# Trace Preserving Homomorphisms on SL(2,C)

By Norman Purzitsky

**Abstract**  Let G, $\Gamma$ be subgroups of SL(2,C), the group of 2x2 matrices of determinant 1 with complex entries. Let $\boldsymbol{\varphi}$: G$\to$ $\boldsymbol{\Gamma}$ be a homomorphism. We call $\boldsymbol{\varphi}$ a trace preserving homomorphism if tr($\boldsymbol{\varphi}$(g))=tr(g) for all g∈G, where tr(g) is the trace of g. We solve the question of when a trace invariant homomorphism is a conjugation by some A∈ SL(2,R). Moreover, if the group G is finitely presented, this paper determines which traces of the generators and products of the generators determine the group up to conjugation. Incomplete solutions are known from the study of Fuchsian groups. Our theorems in this paper will expand the results in the literature to include Fuchsian Groups with elliptic elements, which have not been considered before. Moreover, they will be applicable to any class of subgroups of SL(2,C). The methods used will be relatively elementary and will indicate how many traces are needed, and the role that any relator equation plays in the parameterization by traces of a class of groups.

## 1. Introduction

Let G, $\Gamma$ be subgroups of **SL(2,C)**, the group of 2x2 matrices of determinant 1 with complex entries. Let $\varphi$: G$\to$ $\Gamma$ be a homomorphism. We call $\varphi$ a **trace preserving homomorphism** if tr($\varphi$(g))=tr(g) for all g∈G, where tr(g) is the trace of g. In the parameterization of Fuchsian Groups, i.e. discrete subgroups of **SL(2,R)**, where R is the set of real numbers, two questions are addressed[5].The first is; when is a trace invariant homomorphism a conjugation by some A∈ **SL(2,R)**. The second question is; if the group G is finitely presented, which traces of the generators and products of the generators determine the group up to conjugation. Partial solutions are given in[5]. Our theorems in this paper will expand the results in the literature to include Fuchsian Groups with elliptic elements, which have not been considered before. Moreover, they will be applicable to any class of subgroups of **SL(2,R)**, as well as classes of subgroups of **SL(2,C)**. The methods used will be relatively elementary and will indicate how many parameters, using only traces, are needed, and the role that any relator equation plays in the parameterization of a class of groups.

## 2. Preliminaries.

If g=$\begin{pmatrix} a & b \\ c & d \end{pmatrix}$ ∈G⊂ **SL(2,C)**, we may associate with g a linear fractional transformation g(z)=$\frac{az+b}{cz+d}$, with ad-bc=1. If tr(g) ≠ 2, then g(z) has two fixed points and is therefore conjugate to $\rho^2 z$, where $\rho \neq 1$. Since the map $\begin{pmatrix} a & b \\ c & d \end{pmatrix} \to \frac{az+b}{cz+d}$ is a homomorphism with kernel $\pm I = \pm \begin{pmatrix} 1 & 0 \\ 0 & 1 \end{pmatrix}$, we have that the matrix g can be conjugated to the form $\begin{pmatrix} \rho & 0 \\ 0 & \rho^{-1} \end{pmatrix}$, where $\rho \neq 1$. We will not from this point on make any distinction between the matrix and the linear fractional transformation associated with the given matrix. So, $\begin{pmatrix} \rho & 0 \\ 0 & \rho^{-1} \end{pmatrix}$ is matrix whose fixed points are 0 and ∞.

We write the group $G=\langle g_1, \ldots, g_n \rangle$ to mean, that if $g \in G$, then $g=W(g_1, \ldots, g_n)$, where, $W=W(g_1, \ldots, g_n)$ is a word in the letters $g_1, \ldots, g_n$. A **free group** on the $g_i$'s is a group where every non trivial $W \neq \pm I$, where I is the identity matrix. If for some nontrivial W, we have $W=\pm I$, then W is a **relator.**

## 2. Homomorphisms and Conjugation.

**Theorem 1**. Let G, Γ be subgroups of **SL(2,C)**, and $\varphi: G \to \Gamma$ be a surjective trace preserving homomorphism. If there exists $g_1, g_2 \in G$ which have no common fixed points, then $\varphi(g)= AgA^{-1}$ for all $g \in G$ and some $A \in $ **SL(2, C)**.

**Proof.** We first assume that $g_1$ has exactly two fixed points, so that $tr(g_1) \neq 2$. Since $\varphi$ is trace preserving, then $\varphi(g_1)$ has exactly two fixed points. We conjugate G and Γ, respectively, so that after the two conjugations we may assume

$g_1 = \varphi(g_1) = \begin{pmatrix} \rho & 0 \\ 0 & \rho^{-1} \end{pmatrix}$, where $\rho \neq 1$. Write $g_2 = \begin{pmatrix} a & b \\ c & d \end{pmatrix}$ and

$\varphi(g_2) = \begin{pmatrix} \alpha & \beta \\ \gamma & \delta \end{pmatrix}$. Since $a+d = \alpha + \delta$, and $\rho a + d/\rho = \rho \alpha + \delta/\rho$, we have $a = \alpha$ and $d = \delta$. Moreover, for any pair $g_1, g_2 \in $ **SL(2, C)**, $g_1, g_2$ have a common fixed point if and only if $2 - tr(g_1 g_2 g_1^{-1} g_2^{-1}) = bc(\rho - 1/\rho) = 0$. Since $\varphi$ is trace invariant and $g_1, g_2$ do not share a fixed point, then $\varphi(g_1), \varphi(g_2)$ do not share a fixed point. Therefore, $bc = \beta\gamma$ are different from 0. We may now conjugate G, Γ, respectively, by a transformation, which keeps 0 and ∞ fixed, so that we may, also, assume $c = \gamma = 1$ and $\beta = b$.

Let $g \in G$ with $g = \begin{pmatrix} w & x \\ y & z \end{pmatrix}$ and $\varphi(g) = \begin{pmatrix} w' & x' \\ y' & z' \end{pmatrix}$. By the same calculations above for $g_2$ and $\varphi(g_2)$ we see that $w = w'$ and $z = z'$. Since $\varphi(g)$ preserves traces and $\varphi(g_2) = g_2 = \begin{pmatrix} a & b \\ 1 & d \end{pmatrix}$, we calculate :

$aw' + by' + x' + dz' = aw + by + x + dz$

$aw' + \rho^2 by' + \frac{x'}{\rho^2} + dz' = aw + \rho^2 by + \frac{x}{\rho^2} + dz$.

We now conclude from these equations that $x' = x$ and $y' = y$, i.e. $\varphi(g) = g$ for all $g \in G$, i.e the original G and Γ are conjugates of each other by $\varphi$.

If both $tr(g_1) = tr(g_2) = 2$, but do not share a fixed point, the both are parabolic transformations. We conjugate so that $g_1 = \begin{pmatrix} 1 & 1 \\ 0 & 1 \end{pmatrix}$ and $g_2 = \begin{pmatrix} 1 & 0 \\ \lambda & 1 \end{pmatrix}$, where $\lambda > 0$. Now $tr(g_1 g_2) > 2$ and neither fixed point is 0 or ∞, therefore we can apply the above calculations.

**Example:** Let $A = \begin{pmatrix} \rho & 0 \\ 0 & \rho^{-1} \end{pmatrix}$ and $B = \begin{pmatrix} \sigma & 1 \\ 0 & \sigma^{-1} \end{pmatrix}$. Set $\varphi(A) = \begin{pmatrix} \rho & 0 \\ 0 & \rho^{-1} \end{pmatrix}$,

$\varphi(B) = \begin{pmatrix} \sigma & 0 \\ 0 & \sigma^{-1} \end{pmatrix}$ and $\varphi(W(A,B)) = W(\varphi(A), \varphi(B))$. Note $\varphi$ preserves traces, but is not an isomorphism.

# 3. Trace Algebra or Homomorphisms

The following four facts will be used.
Let $\varphi: G \to \Gamma$ be a homomorphism. Then:

(1) $tr(\varphi(X))=tr(X)$ if and only if $tr(\varphi(X^{-1}))=tr(X^{-1})$;
(2) $tr(\varphi(XYX^{-1}))=tr(\varphi(X))$:
(3) $tr(\varphi(X))=tr(X)$ if and only $tr(\varphi(YXY^{-1}))=tr(YXY^{-1})$ for some X;
(4) The formula
$$tr(XY)+tr(X^{-1}Y)=tr(X)tr(Y)$$

for traces of matrices in **SL(2,C)** implies that if $\varphi$ preserves the trace of any three matrices in this formula, it preserves the trace in the fourth.

**Proposition 1.** If $tr(\varphi(W))=tr(W)$, $tr(\varphi(A))=tr(A)$, and $tr(\varphi(A^{-1}W))=tr(A^{-1}X)$,
Then: (a) $tr(\varphi(AW))=tr(AW)$;
(b) $tr(\varphi(AW^{-1}))=tr(AW^{-1})$;
(c) $tr(\varphi(A^{-1}W^{-1}))=tr(A^{-1}W^{-1})$.

The equations (a), (b), (c) follow easily from the formula in (4) above.

**Proposition 2.** If $\varphi$ preserves the trace of A, B, AB, W, WB$^{-1}$, AW, and ABW, then $\varphi$ preserves the trace of BAW.

**Proof:** Using (4) above we write:

$tr(B^{-1}A^{-1}W)+tr(ABW)=tr(AB)tr(W)$ implies $\varphi$ preserves the trace of B$^{-1}$A$^{-1}$W;

$tr(A^{-1}WB^{-1})=tr(B^{-1}A^{-1}W)$ implies $\varphi$ preserves the trace of $A^{-1}WB^{-1}$;

$tr(AWB^{-1})+tr(A^{-1}WB^{-1})=tr(A)tr(WB^{-1})$ give us that $\varphi$ preserves the trace of $AWB^{-1}$:

$tr(AWB^{-1})=tr(B^{-1}AW)$ implies $\varphi$ preserves the trace of $B^{-1}AW$.
Finally
$tr(BAW)+tr(B^{-1}AW)=tr(B) tr(AW)$ gives the result that $\varphi$ preserves the trace of BAW.

**Corollary:** If $W=W_1W_2$ in Proposition 2 and $\varphi$ preserves the trace of $W_2ABW_1$, instead of ABW, then $\varphi$ preserves the trace of $W_2BAW_1$.
**Proof:** Conjugation by $W_2$ equates the two situations.

**Remark.** The main point of this proposition is that under the appropriate circumstances the $tr(\varphi(W_2ABW_1))=tr(\varphi(W_2BAW_1))$.

3.

# Trace Preserving Homomorphisms

**Theorem 2.** Let $G=\langle A_1 \ldots A_n \rangle$. Let $\varphi: G \to \Gamma$ be a surjective homomorphism such that $\text{tr}(\varphi(A_i))=\text{tr}(A_i)$ for $i=1,2\ldots n$, and $\text{tr}(\varphi(A_{i_1} \ldots A_{i_k}))=\text{tr}(A_{i_1} \ldots A_{i_k})$ for $k \leq n$, where $i_1 < i_2 \ldots < i_k$ for all choices of $i_1, \ldots, i_k$ and $k \leq n$. Then $\text{tr}(\varphi(X))=\text{tr}(X)$ for all $X \in G$.

**Proof:** For each word $W = A_{i_1}^{\alpha_1} \ldots A_{i_k}^{\alpha_k}$ set $\ell(W) = |\alpha_1| + \ldots + |\alpha_k|$. We will proceed by induction on $\ell(W)$. However, we must first do the case where each letter in W occurs exactly once. In this case $W = A_{i_1} \ldots A_{i_k}$, where the $i_j$ are all distinct integers, but not in any particular order. Applying Proposition 2 as many time as needed, we can conclude $\varphi$ preserves the trace of W if and only if it preserves the trace of the word W' with the $i_j$'s arranged in ascending order. Now by hypothesis $\varphi$ preserves the trace of W', and therefore W.

To proceed by induction we note by hypothesis the conclusion is true for $\ell(W)=1$. Suppose $\text{tr}(\varphi(W))=\text{tr}(W)$ for $\ell(W)<n$. Let $\ell(W)=n$, and one of the letters, which we will denote by A, in the generating set occurs twice in W.

There are two cases. For the first case write W=XAYAZ, where X,Y,Z are words in the $A_i$'s. Now $\varphi$ preserves the trace of W if and only if it preserves the trace of AYAZX. By the induction hypothesis $\varphi$ preserves the trace of AY, AZX, and Y⁻¹ZX. Using the formula in (4) above, we conclude that $\varphi$ preserves the trace of AYAZX, and therefore W.

In the second case we write W=XAYA⁻¹Z. Now $\varphi$ preserves the trace of W if and only if it preserves the trace of AYA⁻¹ZX. Now by induction $\varphi$ preserves the trace of AY and A⁻¹ZX. Therefore, $\varphi$ preserves the trace of both AYA⁻¹ZX and AYX⁻¹Z⁻¹A together, or neither of the two traces is preserved. By the first case $\varphi$ preserves the trace of AYX⁻¹Z⁻¹A.

**Corollary:** All the traces of the matrices in $\Gamma$ are polynomials in the $\text{tr}(\varphi(A_i))$ for $i=1,2\ldots n$, and the $\text{tr}(\varphi(A_{i_1} \ldots A_{i_k}))$ for $k \leq n$, where $i_1 < i_2 \ldots < i_k$ for all choices of $i_1, \ldots, i_k$ and $k \leq n$.

The following example shows the necessity of the full hypothesis up to the case n=3.

**Example1:** Let $A=\begin{pmatrix} 3 & 0 \\ 0 & \frac{1}{3} \end{pmatrix}$, $B=\begin{pmatrix} 1 & -3 \\ 1 & -2 \end{pmatrix}$, $C=\begin{pmatrix} 4 & \beta \\ \gamma & -2 \end{pmatrix}$. If G is any free subgroup of **SL(2,C)**, on the 3 generators $A_1, A_2, A_3$, and $\varphi(A_1)=A$, $\varphi(A_2)=B$, $\varphi(A_3)=C$, note tr(AC), tr(AB) are uniquely determined, but tr(BC)= 4-3$\gamma$+$\beta$+4 is not for any value that we set tr(BC). Now the equation -8-$\gamma\beta$=1 limits our solutions to 2 possibilities.

So set tr(BC)=8, the $\gamma = \mp i\sqrt{3}$ and $\beta = \frac{-9}{\gamma}$ are the two solutions. Only by knowing tr(ABC) will the value of $\gamma$ be determined.

**Example 2:** Let $B_i = \begin{pmatrix} \rho_i & \beta_i \\ 0 & \rho_i^{-1} \end{pmatrix}$, for i=1,2...n, and $A_i$, i=1,2...,n be the free generators of subgroup G of **SL(2,C)**. Set $\varphi(A_i)=B_i$ for i=1,2...,n. Note that tr($B_i$) determines $\rho_i$ and $\rho_i^{-1}$ without knowing which is the (1,1) entry. Moreover, note that the $\beta_i$'s are never part of the trace calculation in any word. Therefore tr($B_i B_j$) completely determine the trace of all words in $\Gamma = \langle B_1 \ldots B_n \rangle$.

This brings us to the other situation.

**Theorem 3..** Let $G = \langle A_1 \ldots A_n \rangle$. Let $\varphi: G \to \Gamma$ be a surjective homomorphism such that tr($\varphi(A_i)$)=tr($A_i$) for i=1,2...n, and $\varphi(A_1)$, $\varphi(A_2)$ do not share a fixed point. Moreover, let tr($\varphi(A_1 A_2)$)=tr($A_1 A_2$), tr($\varphi(A_1 A_i)$)=tr($A_1 A_i$), tr($\varphi(A_2 A_i)$)=tr($A_2 A_i$), and tr($\varphi(A_1 A_2 A_i)$)=tr($A_1 A_2 A_i$) for i=3,...,n. Then $\varphi$ is trace preserving.

**Proof:** We first do the case that $\varphi(A_1)$ has exactly two fixed points. We may conjugate $\Gamma$ so that $\varphi(A_1) = \begin{pmatrix} \rho & 0 \\ 0 & \rho^{-1} \end{pmatrix}$, $\varphi(A_2) = \begin{pmatrix} a & b \\ 1 & d \end{pmatrix}$, and write $\varphi(A_i) = \begin{pmatrix} \alpha & \beta \\ \gamma & \delta \end{pmatrix}$.

Now apply calculations very similar to those of **Theorem 1** to determine uniquely all the entries of $A_i$. This makes $\varphi$ trace preserving.

If $A_1$ has only one fixed point then either $A_2$ or $A_1 A_2$ has two fixed points. Using **Proposition 2**, we see that we may switch the roles of $A_1$ and $A_2$, when $A_2$ has two fixed points. If neither $A_1$ nor $A_2$ have two fixed points, then we use $A_1 A_2$ and $A_2^{-1}$ in the calculations above, noting the hypothesis still holds.

## The Relator.

If a word $W(\varphi(A_1) \ldots \varphi(A_n)) = \pm I$ or is a parabolic transformation, then we have the equation tr($W(\varphi(A_1) \ldots \varphi(A_n))) = \pm 2$. It is possible that from this equation we can explicitly solve for one of the traces. However, since tr($W(\varphi(A_1) \ldots \varphi(A_n))$) is a polynomial in the given traces above, we can apply the Implicit Function Theorem, when necessary, to implicitly solve for one of the given traces in terms of the others. There is the proviso, that one of the partial derivatives in the polynomial tr($(\varphi(A_1) \ldots \varphi(A_n))$) be non zero. Moreover, in this case, if the other variables are allowed to run over open sets in **C** or **R**, then we have local coordinates for the class of groups considered as a manifold in $\mathbf{C}^n$ or $\mathbf{R}^n$ for the appropriate field and n.

# References.

Department of Mathematics
York University
North York, Ontario M3J-1P6
Canada
Email: purzit@mathstat.yorku.ca